\documentclass[11pt,leqno]{article}
\usepackage{amsmath,amstext}

\date{31 July 1998}

\parfillskip=\parindent plus1fil
\renewcommand{\phi}{\varphi}

\title{\sffamily\bfseries\Large Commutants of Analytic Toeplitz Operators\\
on the Bergman Space}

\author{\sffamily\bfseries Sheldon Axler,
\v{Z}eljko \v{C}u\v{c}kovi\'{c}, N. V. Rao}

\newcounter{referencec}
\newcommand{\bibref}[1]{[\ref{#1}]}

\newenvironment{keeptogether}{\pagebreak[0]\samepage}{}

\begin{document}
\maketitle
\thispagestyle{empty}

\begin{list}{}{\setlength{\leftmargin}{\parindent}
\setlength{\rightmargin}{\parindent}}

\item \small
\textit{Abstract}. In this note we
show that if two Toeplitz operators on a Bergman space
commute and the symbol of one
of them is analytic and nonconstant, then the other one is also
analytic.
\end{list} 

Let $\Omega$ be a bounded open domain in the complex plane and let
$dA$ denote area measure on $\Omega$. The Bergman space $L_a^2(\Omega)$ is
the subspace of $L^2(\Omega,dA)$ consisting of the square-integrable
functions that are
analytic on $\Omega$. For a bounded measurable
function $\phi$ on $\Omega$, the Toeplitz operator $T_{\phi}$
with symbol $\phi$ is the operator on $L_a^2(\Omega)$ defined
by\footnote{The first author
was partially supported by the National Science Foundation.}
\footnote{\textit{Mathematics Subject Classification}: 47B35}
\[
T_{\phi}(f)=P(\phi f),
\]
where $P$ is the orthogonal projection of $L^2(\Omega,dA)$ onto $L_a^2(\Omega)$.
A~Toeplitz operator is called
analytic if its symbol is an analytic function on~$\Omega$.
Note that if $\phi$ is a bounded analytic function on $\Omega$, then $T_\phi$
is simply the operator of multiplication by $\phi$ on $L_a^2(\Omega)$.

The general problem that we are interested in is the following:
When two Toeplitz operators commute,
what is the relationship between their symbols? If we were working on the
Hardy space of the circle instead of the Bergman space,
then the following result
would answer this question:

\begin{itemize}
\item
If Hardy space Toeplitz operators $T_{\phi}$ and $T_{\psi}$ commute, then either
both symbols are analytic or both symbols are conjugate analytic or
$a \phi + b \psi$ is constant for some constants $a, b$ not both $0$
(Brown and Halmos \bibref{BH}).
\end{itemize}
More general results concerning which operators, not necessarily Toeplitz,
commute with an analytic Hardy space Toeplitz operator are due to
Thomson (\bibref{Th1} and \bibref{Th2}) and Cowen \bibref{Cow}.

On the Bergman space, the situation is more complicated. The Brown-Halmos result
mentioned above fails. For example, if $\Omega$ is the unit disk, then any
two Toeplitz operators whose symbols are radial functions commute (proof:
an easy calculation shows that every Toeplitz operator
with radial symbol has a diagonal matrix with respect to the usual
orthonormal basis; any two diagonal matrices commute).

Despite the difficulty of the general problem, we are
encouraged by the partial results known when $\Omega$ is the unit disk.
If $\Omega$ is the unit disk and
$T_\phi$ and $T_\psi$ commute, then the following hold:

\begin{itemize}
\item
If  $\phi=z^n$,
then $\psi$ is analytic (\v{C}u\v ckovi\'{c} \bibref{Cuc}).

\item
If $\phi$ and $\psi$ are both harmonic, then either both symbols are analytic or
both symbols are conjugate analytic
or $a \phi + b \psi$ is constant for some constants $a, b$ not both $0$
(Axler and \v{C}u\v{c}kovi\'{c} \bibref{AC}).

\item
If  $\phi$ is a radial function, then $\psi$
is radial (\v{C}u\v ckovi\'{c} and Rao \bibref{CR}). 

\item
If $\phi=z^m\bar{z}^n$, then
$\psi(re^{i\theta}) = \sum_{j=-\infty}^\infty \psi_j(r) e^{ij\theta}$, where
$\{\psi_j\}$ are the functions (depending upon $m,n$) described
by \v{C}u\v{c}kovi\'{c} and Rao \bibref{CR}.
\end{itemize}

In this note we extend \v{C}u\v ckovi\'{c}'s first result above by replacing
the disk with an arbitrary bounded domain and (more importantly) by replacing
$z^n$ with an arbitrary bounded analytic function. Here is our result:

\bigskip

\noindent
\textbf{Theorem:} \textsl{If $\phi$ is a nonconstant bounded analytic function on
$\Omega$ and $\psi$ is a bounded measurable function on $\Omega$ such
that $T_\phi$ and $T_\psi$ commute, then $\psi$ is analytic.}

\bigskip

Our proof depends on the following approximation theorem:  

\begin{itemize}
\item
Let $\phi$ be a nonconstant bounded analytic function on $\Omega$.
Then the norm closed subalgebra of $L^{\infty}(\Omega,dA)$ generated by
$\bar{\phi}$ and the bounded analytic functions on $\Omega$
contains $C(\bar{\Omega})$ (Bishop \bibref{Bis}).
\end{itemize}
There is a large literature of related approximation theorems;
see, for example,
\v{C}irka \bibref{Cir}, Axler and Shields \bibref{AS}, Izzo \bibref{Izz}.

\medskip

\textsc{Proof of Theorem:}
Suppose $\phi$ is a nonconstant bounded analytic function on
$\Omega$ and $\psi$ is a bounded measurable function on $\Omega$ such
that $T_{\phi}T_{\psi}=T_{\psi}T_{\phi}$.

Write $\psi=f + u$ with $f\in L_a^2(\Omega)$ and
$u \in L^2(\Omega) \ominus L_a^2(\Omega)$. If $n$ is a nonnegative integer, then
\begin{align*}
T_{\phi^n}T_{\psi}(1) &= \phi^n P(f + u) \\
&= \phi^n f
\end{align*}
and
\begin{align*}
T_{\psi}T_{\phi^n}(1) &= P(f \phi^n + u \phi^n) \\
&= f \phi^n + P(u\phi^n).
\end{align*}
Our hypothesis implies that $T_{\phi^n}T_{\psi} = T_{\psi}T_{\phi^n}$, and
thus the equations above imply that $P(u \phi^n)=0$. Hence if $h \in L_a^2(\Omega)$
we have
\begin{align*}
0 &= \langle h, u \phi^n \rangle \\
&= \int_\Omega \bar{u} h  \overline{\phi^n}\,dA.
\end{align*}
Because the equation above holds for every bounded analytic function $h$ on
$\Omega$ and every nonnegative integer $n$, Bishop's result quoted above
implies that
\[
\int_\Omega \bar{u} w\,dA = 0
\]
for every $w \in C(\bar{\Omega})$. But $C(\bar{\Omega})$ is
dense in $L^2(\Omega,dA)$, and so this implies that $u=0$.
Thus $\psi=f$ and hence $\psi$ is analytic,
completing the proof.

\section*{Open Problems}
\begin{itemize}
\item
If an operator $S$ in the algebra generated
by the Toeplitz operators commutes with a nonconstant analytic Toeplitz operator,
then is $S$ itself Toeplitz and hence (by our result) analytic?

\item
Suppose $\Omega$ is the unit disk and $\phi$ is a bounded harmonic function
on the disk that is
neither analytic nor conjugate analytic. If $\psi$ is a bounded measurable
function on the disk such that $T_\phi$ and $T_\psi$ commute,
must $\psi$ be of the form $a\phi+b$ for some constants $a,b$?
This question would have a negative answer if the disk were replaced by
an annulus centered at the origin because $T_{\log |z|}$ commutes with
every Toeplitz operator with radial symbol.

\item
What is the situation on Bergman spaces in higher dimensions?
\end{itemize}

\section*{References}

\begin{list}{\arabic{referencec}.\hfill}{\usecounter{referencec} 
\setlength{\topsep}{12pt plus 3pt minus 3pt}
\setlength{\partopsep}{0pt}
\setlength{\labelwidth}{1.5\parindent}
\setlength{\labelsep}{0pt}
\setlength{\leftmargin}{1.5\parindent}
\setlength{\parsep}{1ex} }

\raggedright

\item \label{AC}
Sheldon Axler and \v{Z}eljko \v{C}u\v ckovi\'{c},
Commuting Toeplitz Operators with harmonic symbols,
\textsl{Integral Equations Operator Theory} 14 (1991), 1--12.

\item \label{AS}
Sheldon Axler and Allen Shields,
Algebras generated by analytic and harmonic functions,
\textsl{Indiana Univ. Math. J.} 36 (1987), 631--638.

\item \label{Bis}
Christopher J. Bishop,
Approximating continuous functions by holomorphic and harmonic functions,
\textsl{Trans. Amer. Math. Soc.} 311 (1989), 781--811.

\item \label{BH}
Arlen Brown and P. R. Halmos,
Algebraic properties of Toeplitz operators,
\textsl{J. Reine Angew. Math.} 213 (1964), 89--102.

\item \label{Cir}
E. M. \v{C}irka,
Approximation by holomorphic functions on smooth manifolds in $\mathbf{C}^n$,
\textsl{Mat. Sb.} 78 (1969), 101--123. 

\item \label{Cow}
Carl C. Cowen,
The commutant of an analytic Toeplitz operator,
\textsl{Trans. Amer. Math. Soc.} 239 (1978), 1--31.

\item \label{Cuc}
\v{Z}eljko \v{C}u\v ckovi\'{c},
Commutants of Toeplitz operators on the Bergman space,
\textsl{Pacific J. Math.} 162 (1994), 277--285.

\item \label{CR}
\v{Z}eljko \v{C}u\v{c}kovi\'{c}  and N. V. Rao,
Mellin Transform, monomial symbols and commuting Toeplitz operators,
to appear in \textsl{J. Functional Analysis}.

\item \label{Izz}
Alexander J. Izzo,
Uniform algebras generated by holomorphic and pluriharmonic functions,
\textsl{Trans. Amer. Math. Soc.} 339 (1993), 835--847.

\item \label{Th1}
James E. Thomson,
The commutant of a class of analytic Toeplitz operators,
\textsl{Amer. J. Math.} 99 (1977), 522--529.

\item \label{Th2}
James Thomson,
The commutant of a class of analytic Toeplitz operators II,
\textsl{Indiana Univ. Math. J.} 25 (1976), 793--800.

\end{list}

\bigskip

\bigskip

\pagebreak[2]

\begin{keeptogether}
\noindent
{\scshape
Sheldon Axler \\
Department of Mathematics \\
San Francisco State University \\
San Francisco, CA 94132 USA}
\end{keeptogether}

\bigskip

\begin{keeptogether}
\noindent
{\scshape
\v{Z}eljko \v{C}u\v{c}kovi\'{c} and N.V. Rao \\
Department of Mathematics \\
University of Toledo \\
Toledo, OH 43606 USA}
\end{keeptogether}

\begin{tabbing}
\textsl{e-mail}: \=\texttt{axler@sfsu.edu}\\
\>\texttt{zcuckovi@math.utoledo.edu}\\
\>\texttt{rnagise@math.utoledo.edu}
\end{tabbing}

\begin{tabbing}
\textsl{www home pages}:\\
\quad \=\texttt{http://math.sfsu.edu/axler}\\
\>\texttt{http://www.math.utoledo.edu/faculty\_pages/zcuckovic.html}\\
\>\texttt{http://www.math.utoledo.edu/faculty\_pages/rnagisetty.html}
\end{tabbing}

\end{document}